\documentclass{article}
\usepackage{amsmath}
\usepackage{amssymb}

\usepackage{theorem}
\numberwithin{equation}{section}

\newtheorem{theorem}{Theorem}[section]

\newtheorem{lemma}[theorem]{Lemma}
\newtheorem{definition}{Definition}[section]

\begin{document}
\title{Global well-posedness and scattering for the radial, defocusing, cubic wave equation with initial data in a critical Besov space}
\date{\today}
\author{Benjamin Dodson}
\maketitle

\noindent \textbf{Abstract:} In this paper we prove that the cubic wave equation is globally well - posed and scattering for radial initial data lying in $B_{1,1}^{2} \times B_{1,1}^{1}$. This space of functions is a scale invariant subspace of $\dot{H}^{1/2} \times \dot{H}^{-1/2}$.

\section{Introduction}
The cubic nonlinear wave equation in three dimensions,

\begin{equation}\label{1.1}
u_{tt} - \Delta u = -u^{3} = F(u), \hspace{5mm} u(0,x) = u_{0}, \hspace{5mm} u_{t}(0,x) = u_{1}, \hspace{5mm} x \in \mathbf{R}^{3},
\end{equation}
has been a topic of recent interest in the study of dispersive partial differential equations. A solution to $(\ref{1.1})$ has the Hamiltonian

\begin{equation}\label{1.2}
E(u(t)) = \frac{1}{2} \int |\nabla u(t,x)|^{2} dx + \frac{1}{2} \int u_{t}(t,x)^{2} dx + \frac{1}{4} \int u(t,x)^{4} dx.
\end{equation}
A solution to $(\ref{1.1})$ also obeys the scaling symmetry that if $u(t,x)$ solves $(\ref{1.1})$, then for any $\lambda > 0$,

\begin{equation}\label{1.3}
\lambda u(\lambda t, \lambda x)
\end{equation}
solves $(\ref{1.1})$ with initial data $(\lambda u_{0}(\lambda x), \lambda^{2} u_{1}(\lambda x))$. In three dimensions this problem is called $\dot{H}^{1/2}$ - critical because the symmetry $(\ref{1.3})$ preserves the $\dot{H}^{1/2}(\mathbf{R}^{3}) \times \dot{H}^{-1/2}(\mathbf{R}^{3})$ norm of the initial data.\vspace{5mm}

\noindent Using the arguments found in \cite{CCT} one can show that the initial value problem $(\ref{1.1})$ fails to be even locally well - posed for data lying in spaces less regular than $\dot{H}^{1/2} \times \dot{H}^{-1/2}$, that is, any space $\dot{H}^{s} \times \dot{H}^{s - 1}$, $s < \frac{1}{2}$.

\begin{definition}[Locally well - posed]\label{d1.1}
The initial value problem $(\ref{1.1})$ is said to be locally well - posed on an open interval $0 \in I \subset \mathbf{R}$ in $\dot{H}^{s} \times \dot{H}^{s - 1}$ if

\begin{enumerate}
\item A unique solution $u \in L_{t}^{\infty} \dot{H}^{s}(I \times \mathbf{R}^{3}) \cap L_{t, loc}^{4} L_{x}^{4}(I \times \mathbf{R}^{3})$, $u_{t} \in L_{t}^{\infty} \dot{H}^{s - 1}(I \times \mathbf{R}^{3})$ exists,

\item $u$ is continuous in time, $u \in C(I ; \dot{H}^{s}(\mathbf{R}^{3}))$, $u_{t} \in C(I ; \dot{H}^{s - 1}(\mathbf{R}^{3}))$,

\item $u$ depends continuously on the initial data. That is, for any compact $J \subset I$, if $\| u_{0} - u_{0}^{\ast} \|_{\dot{H}^{s}} < \epsilon$ and $\| u_{1} - u_{1}^{\ast}\|_{\dot{H}^{s - 1}} < \epsilon$, for some $\epsilon(J) > 0$ sufficiently small, then 

\begin{equation}\label{1.4}
\| u^{\ast} - u \|_{L_{t,x}^{4}(J \times \mathbf{R}^{3})} + \| u^{\ast} - u \|_{L_{t}^{\infty} \dot{H}^{s} (J \times \mathbf{R}^{3})} + \| u_{t}^{\ast} - u_{t} \|_{L_{t}^{\infty} \dot{H}^{s - 1} (J \times \mathbf{R}^{3})} \leq C(\epsilon),
\end{equation}
where $C(\epsilon)$ is a continuous function of $\epsilon$, $C(0) = 0$. $u$ is the unique solution with initial data $(u_{0}, u_{1})$ and $u^{\ast}$ is the solution with initial data $(u_{0}^{\ast}, u_{1}^{\ast})$.
\end{enumerate}
\end{definition}

\begin{lemma}\label{l1.2}
$(\ref{1.1})$ is locally well - posed in $\dot{H}^{s} \times \dot{H}^{s - 1}$ for any $s \geq \frac{1}{2}$.
\end{lemma}
\emph{Proof:} See \cite{LS}. $\Box$\vspace{5mm}

\noindent Clearly the results of \cite{CCT} and \cite{LS} completely work out the local theory of the initial value problem $(\ref{1.1})$.\vspace{5mm}

\noindent For the global theory, an additional obstacle is the lack of a conserved Hamiltonian (like $(\ref{1.2})$) at the critical regularity. Indeed \cite{Gril} showed that the energy - critical problem obtained either by changing $-u^{3}$ to $-u^{5}$ in $(\ref{1.1})$ or by analyzing $(\ref{1.1})$ in dimension four is globally well - posed. The proof argues by showing that the conserved energy cannot concentrate at the tip of a light cone.\vspace{5mm}

\noindent \textbf{Remark:} The focusing case (replace $-u^{3}$ in $(\ref{1.1})$ by $u^{5}$) is considerably more complicated. Focusing problems are not addressed at all in this paper, and so the interested reader will simply be referred to \cite{Kenig} and the references therein.\vspace{5mm}

\noindent However, there is no known conserved quantity that controls the $\dot{H}^{1/2} \times \dot{H}^{-1/2}$ norm of a solution to $(\ref{1.1})$. For the radial version of $(\ref{1.1})$,

\begin{equation}\label{1.5}
u_{tt} - u_{rr} - \frac{2}{r} u_{r} + u^{3} = 0, \hspace{5mm} u(0,x) = u_{0}, \hspace{5mm} u_{t}(0,x) = u_{1},
\end{equation}
this is in fact the obstacle to proving that $(\ref{1.1})$ is globally well - posed and scattering in $\dot{H}^{1/2} \times \dot{H}^{-1/2}$.

\begin{theorem}\label{t1.3}
Suppose $u_{0} \in \dot{H}^{1/2}(\mathbf{R}^{3})$, $u_{1} \in \dot{H}^{-1/2}(\mathbf{R}^{3})$ are radial functions, and $u$ solves $(\ref{1.1})$ on a maximal interval $0 \in I \subset \mathbf{R}$ with

\begin{equation}\label{1.6}
\sup_{t \in I} \| u(t) \|_{\dot{H}^{1/2}(\mathbf{R}^{3})} + \| u_{t}(t) \|_{\dot{H}^{-1/2}(\mathbf{R}^{3})} < \infty.
\end{equation}
Then $I = \mathbf{R}$ and the solution $u$ scatters both forward and backward in time.
\end{theorem}
\emph{Proof:} See \cite{DL}. $\Box$

\begin{definition}[Scattering]\label{d1.4}
A solution to $(\ref{1.1})$ is said to scatter forward in time if there exist some $u_{0}^{+} \in \dot{H}^{1/2}$, $u_{1}^{+} \in \dot{H}^{-1/2}$ such that

\begin{equation}\label{1.7}
\lim_{t \rightarrow \infty} \| u(t) - S(t)(u_{0}^{+}, u_{1}^{+}) \|_{\dot{H}^{1/2}(\mathbf{R}^{3})} + \| u_{t}(t) - \partial_{t} S(t)(u_{0}^{+}, u_{1}^{+}) \|_{\dot{H}^{-1/2}(\mathbf{R}^{3})} = 0,
\end{equation}
where $u(t) = S(t)(u_{0}, u_{1})$ is the solution to the linear wave equation

\begin{equation}\label{1.8}
u_{tt} - \Delta u = 0, \hspace{5mm} u(0,x) = u_{0}, \hspace{5mm} u_{t}(0,x) = u_{1}.
\end{equation}
A solution to $(\ref{1.1})$ is said to scatter backward in time if there exist $u_{0}^{-} \in \dot{H}^{1/2}$, $u_{1}^{-} \in \dot{H}^{-1/2}$ such that

\begin{equation}\label{1.9}
\lim_{t \rightarrow -\infty} \| u(t) - S(t)(u_{0}^{-}, u_{1}^{-}) \|_{\dot{H}^{1/2}(\mathbf{R}^{3})} + \| u_{t}(t) - \partial_{t} S(t)(u_{0}^{-}, u_{1}^{-}) \|_{\dot{H}^{-1/2}(\mathbf{R}^{3})} = 0.
\end{equation}
\end{definition}
To compensate for the lack of a quantity that controls the $\dot{H}^{1/2} \times \dot{H}^{-1/2}$ norm, we will use the fact that the cubic exponent in $(\ref{1.1})$ is the conformal exponent ($\frac{n + 3}{n - 1}$) in three dimensions. Working in hyperbolic coordinates \cite{Tataru} proved weighted Strichartz estimates that extended previous results of \cite{GLS}.\vspace{5mm}

\noindent Also working in hyperbolic coordinates \cite{Shen} was able to prove a scattering result for data lying in a weighted energy space. To do this \cite{Shen} used a Morawetz estimate in hyperbolic space. \cite{D} combined the result of \cite{Shen} with the I - method, proving

\begin{theorem}\label{t1.5}
Suppose there exists a positive constant $\epsilon > 0$ such that

\begin{equation}\label{1.10}
\| u_{0} \|_{\dot{H}^{1/2 + \epsilon}(\mathbf{R}^{3})} + \| |x|^{2 \epsilon} u_{0} \|_{\dot{H}^{1/2 + \epsilon}(\mathbf{R}^{3})} \leq A < \infty,
\end{equation}
and

\begin{equation}\label{1.11}
\| u_{1} \|_{\dot{H}^{-1/2 + \epsilon}(\mathbf{R}^{3})} + \| |x|^{2 \epsilon} u_{1} \|_{\dot{H}^{-1/2 + \epsilon}(\mathbf{R}^{3})} \leq A < \infty.
\end{equation}
Then $(\ref{1.1})$ has a global solution and there exists some $C(A, \epsilon) < \infty$ such that

\begin{equation}\label{1.12}
\int_{\mathbf{R}} \int (u(t,x))^{4} dx dt \leq C(A, \epsilon),
\end{equation}
which proves that $u$ scatters both forward and backward in time.
\end{theorem}

\noindent \textbf{Remark:} A straightforward application of the Strichartz estimates of \cite{GV} and \cite{Stri} shows that

\begin{equation}\label{1.13}
\| u \|_{L_{t,x}^{4}(\mathbf{R} \times \mathbf{R}^{3})} < \infty
\end{equation}
is equivalent to scattering.\vspace{5mm}

\noindent Notice that conditions $(\ref{1.10})$ and $(\ref{1.11})$ fall just short of lying in the critical Sobolev space $\dot{H}^{1/2} \times \dot{H}^{-1/2}$, and are not invariant under the scaling $(\ref{1.3})$. In this paper we will study the radial, nonlinear wave equation in three dimensions,

\begin{equation}\label{1.14}
u_{tt} - u_{rr} - \frac{2}{r} u_{r} + u^{3} = 0, \hspace{5mm} u_{0} \in B_{1,1}^{2}, \hspace{5mm} u_{1} \in B_{1,1}^{1}.
\end{equation}
The Besov spaces $B_{q, r}^{s}$ will be defined in the next section. By the Sobolev embedding theorem, this space is a subspace of $\dot{H}^{1/2} \times \dot{H}^{-1/2}$, and the norm is invariant under $(\ref{1.3})$.\vspace{5mm}

\noindent The author believes that this is the first result in which large data scattering was proved for initial data in a scale - invariant space for which the norm was not controlled by a conserved quantity.

\begin{theorem}\label{t1.6}
The initial value problem $(\ref{1.1})$ is globally well - posed and scattering for $u_{0} \in B_{1,1}^{2}(\mathbf{R}^{3})$, radial, and $u_{1} \in B_{1,1}^{1}(\mathbf{R}^{3})$, radial. Moreover,

\begin{equation}\label{1.15}
\| u \|_{L_{t,x}^{4}(\mathbf{R} \times \mathbf{R}^{3})} \leq C(\| u_{0} \|_{B_{1,1}^{2}}, \| u_{1} \|_{B_{1,1}^{1}}).
\end{equation}
\end{theorem}
The proof of this theorem utilizes the fact that the free solution with such initial data is only singular at the origin $t = 0$, $x = 0$. Thus, using a Gronwall - type inequality, the local solution to $(\ref{1.1})$ can be extended to a global solution that is the sum of a solution to the free wave equation combined with a finite energy term. A Morawetz estimate in hyperbolic coordinates then proves scattering.\vspace{5mm}

\noindent The proof of theorem $\ref{t1.6}$ will occupy the remainder of this paper. In section two we will begin by defining the Besov spaces and recalling basic Strichartz estimates. Then in section three the local theory of $(\ref{1.1})$ will be discussed. Global well - posedness will then be proved in section four. In section five we will switch to hyperbolic coordinates to prove scattering. Finally in section six we will use a profile decomposition to show that the bounds obtained for any $u_{0} \in B_{1,1}^{2}$, $u_{1} \in B_{1,1}^{1}$ depend only on size.

\section{Besov spaces and linear estimates}
In this section we present some harmonic analysis estimates that will be used in this paper. None of these results are new.

\begin{theorem}[Hardy - Littlewood - Sobolev inequality]
For any $0 < s < 1$, if $\frac{1}{q} = \frac{1}{p} + s - 1$, then

\begin{equation}\label{2.1}
\| \frac{1}{|t|^{s}} \ast F(t) \|_{L^{q}(\mathbf{R})} \lesssim_{s} \| F \|_{L^{q}(\mathbf{R})}.
\end{equation}
\end{theorem}

\begin{definition}[Littlewood - Paley decomposition]
Let $\phi \in C_{0}^{\infty}(\mathbf{R}^{3})$ be a function supported on $|x| \leq 2$ and $\phi(x) = 1$ for $|x| \leq 1$. Then for any $j \in \mathbf{Z}$ let

\begin{equation}\label{2.2}
P_{j} f = \mathcal F^{-1} (\phi(2^{-j} \xi) \hat{f}(\xi)),
\end{equation}
where

\begin{equation}\label{2.3}
\hat{f}(\xi) = (2 \pi)^{-d/2} \int e^{-ix \cdot \xi} f(x) dx,
\end{equation}
and

\begin{equation}\label{2.4}
\mathcal F^{-1} g = (2 \pi)^{-d/2} \int e^{ix \cdot \xi} g(\xi) d\xi.
\end{equation}
Then for any Schwartz function $f$,

\begin{equation}\label{2.5}
f = \sum_{j \in \mathbf{Z}} P_{j} f.
\end{equation}
\end{definition}
Let $K_{j}(x)$ be the kernel of the Littlewood - Paley multiplier $P_{j}$. Then by direct computation, for any $N$,

\begin{equation}\label{2.6}
|K_{j}(x)| \lesssim_{d, N} \frac{2^{jd}}{(1 + 2^{j} |x|)^{N}}.
\end{equation}
Since $K_{j}$ has an $L^{1}$ norm that is uniformly bounded in $j$, for any $1 \leq p \leq \infty$,

\begin{equation}\label{2.7}
\| P_{j} f \|_{L^{p}(\mathbf{R}^{d})} \lesssim_{d} \| f \|_{L^{p}(\mathbf{R}^{d})}.
\end{equation}
A direct computation also gives Bernstein's inequality

\begin{equation}\label{2.8}
\| P_{j} f \|_{L^{p}(\mathbf{R}^{d})} \lesssim_{d} 2^{-j} \| \nabla f \|_{L^{p}(\mathbf{R}^{d})},
\end{equation}
along with the Sobolev embedding estimate, for $1 \leq p \leq q \leq \infty$,

\begin{equation}\label{2.9}
\| P_{j} f \|_{L^{q}(\mathbf{R}^{d})} \lesssim_{d} 2^{jd(\frac{1}{p} - \frac{1}{q})} \| f \|_{L^{p}(\mathbf{R}^{d})}.
\end{equation}
The Littlewood - Paley decomposition is foundational to the definition of Besov spaces.

\begin{definition}[Besov spaces]\label{d2.0}
Suppose $1 \leq p \leq \infty$, $1 \leq r \leq \infty$, and $s \in \mathbf{R}$. Then

\begin{equation}\label{2.10}
\| f \|_{B_{r, p}^{s}(\mathbf{R}^{d})} = (\sum_{j \in \mathbf{Z}} 2^{jsr} \| P_{j} f \|_{L^{p}(\mathbf{R}^{d})}^{r})^{1/r}
\end{equation}
The Besov space $B_{r,p}^{s}$ is then the completion of the Schwartz space under this norm. $B_{r,p}^{s}$ is a Banach space under this topology.
\end{definition}
The Besov spaces are well - behaved with respect to multiplying by smooth cutoff functions.

\begin{lemma}\label{l2.4}
Suppose $\chi(x) \in C_{0}^{\infty}(\mathbf{R}^{3})$. Then

\begin{equation}\label{2.11}
\| \chi(x) u \|_{B_{1, 2}^{1/2}(\mathbf{R}^{3})} \lesssim \| u \|_{B_{1,2}^{1/2}(\mathbf{R}^{3})},
\end{equation}
and

\begin{equation}\label{2.12}
\| \chi(x) u \|_{B_{1,2}^{-1/2}(\mathbf{R}^{3})} \lesssim \| u \|_{B_{1,2}^{-1/2}(\mathbf{R}^{3})}.
\end{equation}
Also if $\chi(x) = 1$ on $|x| \leq 1$ then if $u_{0} \in B_{1,1}^{2}$ and $u_{1} \in B_{1,1}^{1}$, then

\begin{equation}\label{2.13}
\lim_{R \rightarrow \infty} \| (1 - \chi(\frac{x}{R})) u_{0} \|_{B_{1,2}^{1/2}(\mathbf{R}^{3})} + \| (1 - \chi(\frac{x}{R})) u_{1} \|_{B_{1,2}^{-1/2}(\mathbf{R}^{3})} = 0.
\end{equation}
\end{lemma}
\emph{Proof:} Splitting

\begin{equation}\label{2.13.1}
P_{j}(\chi f) = \chi (P_{j} f) + [P_{j}, \chi] f,
\end{equation}
since by Holder's inequality,

\begin{equation}\label{2.13.2}
\sum_{j} 2^{j/2} \| \chi (P_{j} f) \|_{L^{2}} \lesssim \sum_{j} 2^{j/2} \| P_{j} f \|_{L^{2}},
\end{equation}
it remains to compute

\begin{equation}\label{2.13.3}
\sum_{j} 2^{j/2} \| [P_{j}, \chi] f \|_{L^{2}}.
\end{equation}
By $(\ref{2.6})$ and the fundamental theorem of calculus,

\begin{equation}\label{2.13.4}
\int K_{j}(x - y) [\chi(y) f(y) - \chi(x) f(y)] \lesssim \sum_{k} \int |K_{j}(x - y)| |x - y| |P_{k} f(y)| dy,
\end{equation}
and therefore by Bernstein's inequality and $(\ref{2.6})$,

\begin{equation}\label{2.13.5}
\sum_{j \geq 0} 2^{j/2} \| P_{j}(\chi P_{\geq 0} f) \|_{L^{2}} \lesssim \sum_{j \geq 0} 2^{-j/2} \sum_{k \geq 0} 2^{-k/2} \| P_{k} f \|_{\dot{H}^{1/2}} + \sum_{j} 2^{j/2} \| P_{j} f \|_{\dot{H}^{1/2}} \lesssim \| f \|_{B_{1,2}^{1/2}}.
\end{equation}
Also by Bernstein's inequality, the Sobolev embedding theorem, and H{\"o}lder's inequality,

\begin{equation}\label{2.13.6}
\sum_{j \geq 0} 2^{j/2} \| P_{j} (\chi P_{\leq 0} f) \|_{L^{2}} \lesssim \| \nabla (\chi (P_{\leq 0} f)) \|_{L^{2}} \lesssim \| f \|_{B_{1,2}^{1/2}}.
\end{equation}
Similarly,

\begin{equation}\label{2.13.7}
\| P_{\leq 0} (\chi (P_{\leq 0} f)) \|_{B_{1,2}^{1/2}} \lesssim \| P_{\leq 0} f \|_{L^{6}} \lesssim \| f \|_{B_{1,2}^{1/2}}.
\end{equation}
Finally,

\begin{equation}\label{2.13.8}
\| P_{\leq 0} (\chi (P_{\geq 0} f)) \|_{B_{1,2}^{1/2}} \lesssim \sum_{k \geq 0} 2^{-k/2} \| P_{k} f \|_{\dot{H}^{1/2}} \lesssim \| f \|_{B_{1,2}^{1/2}}.
\end{equation}
Combining $(\ref{2.13.5})$ - $(\ref{2.13.8})$, we have proved

\begin{equation}\label{2.13.9}
\| \chi f \|_{B_{1,2}^{1/2}(\mathbf{R}^{3})} \lesssim \| f \|_{B_{1,2}^{1/2}(\mathbf{R}^{3})}.
\end{equation}
Also observe that $(\ref{2.13.5})$ - $(\ref{2.13.8})$ also imply that

\begin{equation}\label{2.13.10}
\| \chi f \|_{B_{\infty, 2}^{1/2}} \lesssim \| f \|_{B_{\infty, 2}^{1/2}},
\end{equation}
and therefore by duality

\begin{equation}\label{2.13.11}
\| \chi g \|_{B_{1,2}^{-1/2}} \lesssim \| g \|_{B_{1,2}^{-1/2}}.
\end{equation}
To prove $(\ref{2.13})$ observe that $B_{1,2}^{1/2} \times B_{1,2}^{-1/2}$ is invariant under the scaling $(\ref{1.3})$, that is,

\begin{equation}\label{2.13.12}
\aligned
\| (1 - \chi(\frac{x}{R})) u_{0} \|_{B_{1,2}^{1/2}} + \| (1 - \chi(\frac{x}{R})) u_{1} \|_{B_{1,2}^{-1/2}} \\ = R \| (1 - \chi(x)) u_{0}(Rx) \|_{B_{1,2}^{1/2}} + R^{2} \| (1 - \chi(x)) u_{1}(Rx) \|_{B_{1,2}^{-1/2}}.
\endaligned
\end{equation}
The dominated convergence theorem, $(\ref{2.13.6})$, and $(\ref{2.13.7})$ imply that

\begin{equation}\label{2.13.13}
\lim_{R \rightarrow \infty} R \| (1 - \chi(x)) P_{\leq 0} (u_{0}(Rx)) \|_{B_{1,2}^{1/2}} + R^{2} \| (1 - \chi(x)) P_{\leq 0} (u_{1}(Rx)) \|_{B_{1,2}^{-1/2}} = 0.
\end{equation}
Meanwhile, $(\ref{2.13.5})$, $(\ref{2.13.8})$, and the dominated convergence theorem imply that

\begin{equation}\label{2.13.13}
\aligned
\lim_{R \rightarrow \infty} R \| (1 - \chi(x)) P_{\geq 0} (u_{0}(Rx)) \|_{B_{1,2}^{1/2}} + R^{2} \| (1 - \chi(x)) P_{\geq 0} (u_{1}(Rx)) \|_{B_{1,2}^{-1/2}} \\
= \lim_{R \rightarrow \infty} R \sum_{j \geq 0} 2^{j/2} \| (1 - \chi(x)) P_{j}(u_{0}(Rx)) \|_{L^{2}} + R^{2} \sum_{j \geq 0} 2^{-j/2} \| (1 - \chi(x)) P_{j}(u_{1}(Rx)) \|_{L^{2}} = 0.
\endaligned
\end{equation}
$\Box$

\begin{theorem}[Radial Sobolev embedding theorem]\label{t2.3}
For any $j$,

\begin{equation}\label{2.14}
\| |x| P_{j} f \|_{L^{\infty}(\mathbf{R}^{3})} \lesssim \| P_{j} f \|_{\dot{H}^{1/2}(\mathbf{R}^{3})}.
\end{equation}
\end{theorem}
\emph{Proof:} By stationary phase computations, if $f$ is radial then

\begin{equation}\label{2.14.1}
\aligned
f(x) = \int_{0}^{\infty} |\xi|^{2} \hat{f}(|\xi|) \int_{-\pi/2}^{\pi/2} e^{i |x| |\xi| \sin \theta} \cos \theta d\theta d|\xi| \\
= \int_{0}^{\infty} r^{2} \hat{f}(r) \int_{-1}^{1} e^{i |x| r u} du dr = \frac{1}{i |x|} \int_{0}^{\infty} \hat{f}(r) r [e^{i |x| r} - e^{-i |x| r}] dr.
\endaligned
\end{equation}
The theorem then follows by the one dimensional Sobolev embedding theorem $B_{1,2}^{1/2}(\mathbf{R}) \subset L^{\infty}(\mathbf{R})$. $\Box$\vspace{5mm}

\noindent Now observe that the solution to the free wave equation

\begin{equation}\label{2.15}
u_{tt} - \Delta u = 0, \hspace{5mm} u(0,x) = f(x), \hspace{5mm} u_{t}(0,x) = g(x),
\end{equation}
is given by the Fourier multiplier

\begin{equation}\label{2.16}
u(t,x) = \mathcal F^{-1}(\cos(t |\xi|) \hat{f}(\xi) + \frac{\sin(t |\xi|)}{|\xi|} \hat{g}(\xi)) = S(t)(f, g).
\end{equation}

\noindent Then the solution to

\begin{equation}\label{2.18}
u_{tt} - \Delta u = F, \hspace{5mm} u(0,x) = f(x), \hspace{5mm} u_{t}(0,x) = g(x),
\end{equation}
is given by

\begin{equation}\label{2.19}
S(t)(f, g) + \int_{0}^{t} S(t - \tau)(0, F) d\tau.
\end{equation}
\textbf{Remark:} Sometimes, if $u = S(t)(f, g)$ it is convenient to write

\begin{equation}\label{2.20}
(u(t), \partial_{t} u(t)) = S(t)(f, g).
\end{equation}
By standard stationary phase calculations,

\begin{theorem}[Dispersive estimate]\label{t2.2}
\begin{equation}\label{2.20}
\| S(t)(f, g) \|_{L^{\infty}(\mathbf{R}^{3})} \lesssim \frac{1}{t} [\| \nabla^{2} f \|_{L^{1}(\mathbf{R}^{3})} + \| \nabla g \|_{L^{1}(\mathbf{R}^{3})}].
\end{equation}
\end{theorem}
The dispersive estimates can be used to prove Strichartz estimates.

\begin{theorem}\label{t2.1}
Let $I \subset \mathbf{R}$, $t_{0} \in I$, be an interval and let $u$ solve the linear wave equation

\begin{equation}\label{2.21}
u_{tt} - \Delta u = F, \hspace{5mm} u(t_{0}) = u_{0}, u_{t}(t_{0}) = u_{1}.
\end{equation}
Then we have the estimates

\begin{equation}\label{2.22}
\aligned
\| u \|_{L_{t}^{p} L_{x}^{q}(I \times \mathbf{R}^{3})} + \| u \|_{L_{t}^{\infty} \dot{H}^{s}(I \times \mathbf{R}^{3})} + \| u_{t} \|_{L_{t}^{\infty} \dot{H}^{s - 1}(I \times \mathbf{R}^{3})} \\
\lesssim_{p, q, s, \tilde{p}, \tilde{q}} \| u_{0} \|_{\dot{H}^{s}(\mathbf{R}^{3})} + \| u_{1} \|_{\dot{H}^{s - 1}(\mathbf{R}^{3})} + \| F \|_{L_{t}^{\tilde{p}'} L_{x}^{\tilde{q}'}(I \times \mathbf{R}^{3})},
\endaligned
\end{equation}
whenever $s \geq 0$, $2 \leq p, \tilde{p} \leq \infty$, $2 \leq q, \tilde{q} < \infty$, and

\begin{equation}\label{2.23}
\frac{1}{p} + \frac{1}{q} \leq \frac{1}{2}, \hspace{5mm} \frac{1}{\tilde{p}} + \frac{1}{\tilde{q}} \leq \frac{1}{2}.
\end{equation}
\end{theorem}
\emph{Proof:} See for example \cite{Tao}. $\Box$\vspace{5mm}

\noindent \textbf{Remark:} This theorem can easily be combined with the Christ - Kiselev lemma (see \cite{SSog}) and the fact that $|\nabla|$ commutes with the operator $(\partial_{tt} - \Delta)$ to prove many additional estimates.

\begin{lemma}[Perturbation lemma]\label{l2.5}
Let $I \subset \mathbf{R}$ be a time interval. Let $t_{0} \in I$, $(u_{0}, u_{1}) \in \dot{H}^{1/2} \times \dot{H}^{-1/2}$ and some constants $M$, $A$, $A' > 0$. Let $\tilde{u}$ solve the equation

\begin{equation}\label{2.24}
(\partial_{tt} - \Delta) \tilde{u} = F(\tilde{u}) = e,
\end{equation}
on $I \times \mathbf{R}^{3}$, and also suppose $\sup_{t \in I} \| (\tilde{u}(t), \partial_{t} \tilde{u}(t)) \|_{\dot{H}^{1/2} \times \dot{H}^{-1/2}} \leq A$, $\| \tilde{u} \|_{L_{t,x}^{4}(I \times \mathbf{R}^{3})} \leq M$,

\begin{equation}\label{2.25}
\| (u_{0} - \tilde{u}(t_{0}), u_{1} - \partial_{t} \tilde{u}(t_{0})) \|_{\dot{H}^{1/2} \times \dot{H}^{-1/2}} \leq A',
\end{equation}
and

\begin{equation}\label{2.26}
\| e \|_{L_{t,x}^{4/3}(I \times \mathbf{R}^{3})} + \| S(t - t_{0})(u_{0} - \tilde{u}(t_{0}), u_{1} - \partial_{t} \tilde{u}(t_{0})) \|_{L_{t,x}^{4}(I \times \mathbf{R}^{3})} \leq \epsilon.
\end{equation}
Then there exists $\epsilon_{0}(M, A, A')$ such that if $0 < \epsilon < \epsilon_{0}$ then there exists a solution to $(\ref{1.1})$ on $I$ with $(u(t_{0}), \partial_{t} u(t_{0})) = (u_{0}, u_{1})$, $\| u \|_{L_{t,x}^{4}(I \times \mathbf{R}^{3})} \leq C(M, A, A')$, and for all $t \in I$,

\begin{equation}\label{2.27}
\| (u(t), \partial_{t} u(t)) - (\tilde{u}(t), \partial_{t} \tilde{u}(t)) \|_{\dot{H}^{1/2} \times \dot{H}^{-1/2}} \leq C(A, A', M)(A' + \epsilon).
\end{equation}
\end{lemma}
\emph{Proof:} The method of proof is by now fairly well - known. See for example lemma $2.20$ of \cite{KM}. $\Box$

\section{Local theory}
By the dominated convergence theorem, for any $u_{0} \in B_{1,1}^{2}$, $u_{1} \in B_{1,1}^{1}$, and $\delta > 0$ there exists some $j_{0}(\delta) < \infty$ such that

\begin{equation}\label{3.1}
\sum_{j \geq j_{0}} 2^{2j} \| P_{j} u_{0} \|_{L^{1}(\mathbf{R}^{3})} + \sum_{j \geq j_{0}} 2^{j} \| P_{j} u_{1} \|_{L^{1}(\mathbf{R}^{3})} < \delta.
\end{equation}
Then by the rescaling $(\ref{1.3})$ with $\lambda = 2^{-j}$,

\begin{equation}\label{3.2}
\sum_{j \geq 0} 2^{2j} \| P_{j} u_{0} \|_{L^{1}(\mathbf{R}^{3})} + \sum_{j \geq 0} 2^{j} \| P_{j} u_{1} \|_{L^{1}(\mathbf{R}^{3})} < \delta.
\end{equation}

\begin{lemma}\label{l3.1}
Fix $\epsilon_{0} > 0$ small. There exists some $\delta(\epsilon, \| u_{0} \|_{B_{1,1}^{2}}, \| u_{1} \|_{B_{1,1}^{1}}) > 0$ such that

\begin{equation}\label{3.3}
\| u \|_{L_{t,x}^{4}([-\delta, \delta] \times \mathbf{R}^{3})} \lesssim \epsilon,
\end{equation}
and

\begin{equation}\label{3.4}
\| u \|_{L_{t}^{\infty} B_{1,2}^{1/2}([-\delta, \delta] \times \mathbf{R}^{3})} \lesssim 1.
\end{equation}
\end{lemma}
\emph{Proof:} Assume that $(\ref{3.2})$ holds for some $\delta_{1} << \epsilon_{0}$. By the Sobolev embedding theorem and definition $\ref{d2.0}$,

\begin{equation}\label{3.5}
\| S(t)(P_{\leq 0} u_{0}, P_{\leq 0} u_{1}) \|_{L_{x}^{4}(\mathbf{R}^{3})} \lesssim 1,
\end{equation}
while by theorem $\ref{t2.1}$, $(\ref{3.1})$, and $(\ref{3.2})$,

\begin{equation}\label{3.6}
\| S(t)(P_{\geq 0} u_{0}, P_{\geq 0}) \|_{L_{t,x}^{4}(\mathbf{R} \times \mathbf{R}^{3})} \lesssim \delta_{1}.
\end{equation}
Taking $\delta > 0$ sufficiently small, $(\ref{3.5})$ and $(\ref{3.6})$ imply that

\begin{equation}\label{3.7}
\| S(t)(u_{0}, u_{1}) \|_{L_{t,x}^{4}([-\delta, \delta] \times \mathbf{R}^{3})} \lesssim \epsilon_{0}.
\end{equation}
Then by the contraction mapping principle and theorem $\ref{t2.1}$,

\begin{equation}\label{3.8}
\| u \|_{L_{t,x}^{4}([-\delta, \delta] \times \mathbf{R}^{3})} \lesssim \| S(t)(u_{0}, u_{1}) \|_{L_{t,x}^{4}([-\delta, \delta] \times \mathbf{R}^{3})} + \| u \|_{L_{t,x}^{4}([-\delta, \delta] \times \mathbf{R}^{3})}^{3},
\end{equation}
which when $\epsilon_{0} > 0$ is sufficiently small implies

\begin{equation}\label{3.9}
\| u \|_{L_{t,x}^{4}([-\delta, \delta] \times \mathbf{R}^{3})} \lesssim \epsilon_{0}.
\end{equation}
Next observe that by theorem $\ref{t2.1}$ we also have

\begin{equation}\label{3.10}
\| |\nabla|^{1/4} u \|_{L_{t}^{8} L_{x}^{8/3}([-\delta, \delta] \times \mathbf{R}^{3})} + \| |\nabla|^{-1/4} u \|_{L_{t}^{8/3} L_{x}^{8}([-\delta, \delta] \times \mathbf{R}^{3})} \lesssim \epsilon_{0},
\end{equation}
and

\begin{equation}\label{3.11}
\aligned
\| P_{j} u \|_{L_{t,x}^{4}([-\delta, \delta] \times \mathbf{R}^{3})} + \| P_{j} u \|_{L_{t}^{\infty} \dot{H}^{1/2}([-\delta, \delta] \times \mathbf{R}^{3})} \lesssim \| P_{j} u_{0} \|_{\dot{H}^{1/2}(\mathbf{R}^{3})} + \| P_{j} u_{1} \|_{\dot{H}^{-1/2}(\mathbf{R}^{3})} \\ + 2^{-j/2} \sum_{j_{1} \leq j_{2} \leq j - 5} \sum_{j - 3 \leq j_{3} \leq j + 3} \| P_{j_{1}} \|_{L_{t}^{8/3} L_{x}^{8}([-\delta, \delta] \times \mathbf{R}^{3})} \| P_{j_{2}} u \|_{L_{t}^{8/3} L_{x}^{8}([-\delta, \delta] \times \mathbf{R}^{3})} \| P_{j_{3}} u \|_{L_{t,x}^{4}([-\delta, \delta] \times \mathbf{R}^{3})} \\ + 2^{j/4} \sum_{j - 5 \leq j_{1} \leq j_{2} \leq j_{3}} \| P_{j_{1}} u \|_{L_{t,x}^{4}([-\delta, \delta] \times \mathbf{R}^{3})} \| P_{j_{2}} u \|_{L_{t,x}^{4}([-\delta, \delta] \times \mathbf{R}^{3})} \| P_{j_{3}} u \|_{L_{t}^{8} L_{x}^{8/3}([-\delta, \delta] \times \mathbf{R}^{3})},
\endaligned
\end{equation}
so then by $(\ref{3.9})$ and $(\ref{3.10})$,

\begin{equation}\label{3.12}
\aligned
\sum_{j} \| P_{j} u \|_{L_{t,x}^{4}([-\delta, \delta] \times \mathbf{R}^{3})} + \| P_{j} u \|_{L_{t}^{\infty} \dot{H}^{1/2}([-\delta, \delta] \times \mathbf{R}^{3})} \\ \lesssim \sum_{j} \| P_{j} u_{0} \|_{\dot{H}^{1/2}(\mathbf{R}^{3})} + \| P_{j} u_{1} \|_{\dot{H}^{-1/2}(\mathbf{R}^{3})}
+ \epsilon_{0}^{2} \sum_{j}  \| P_{j} u \|_{L_{t,x}^{4}([-\delta, \delta] \times \mathbf{R}^{3})},
\endaligned
\end{equation}
which also implies

\begin{equation}\label{3.13}
\| u \|_{L_{t}^{\infty} B_{1, 2}^{1/2}([-\delta, \delta] \times \mathbf{R}^{3})} \lesssim \| u_{0} \|_{B_{1,1}^{2}(\mathbf{R}^{3})} + \| u_{1} \|_{B_{1,1}^{1}(\mathbf{R}^{3})}.
\end{equation}
$\Box$\vspace{5mm}

\noindent Next suppose $\chi(x)$ is a smooth function that is supported on $|x| \leq 1$ and is equal to one on $|x| \leq \frac{1}{2}$. By lemma $\ref{l2.4}$ there exists some $R(u_{0}, u_{1}, \epsilon)$ such that

\begin{equation}\label{3.14}
\| (1 - \chi(\frac{x}{R})) u_{0} \|_{\dot{H}^{1/2}(\mathbf{R}^{3})} + \| (1 - \chi(\frac{x}{R})) u_{1} \|_{\dot{H}^{-1/2}(\mathbf{R}^{3})} \leq \epsilon.
\end{equation}
\textbf{Remark:} Notice that $R$ depends on $u_{0}$ and $u_{1}$, not just their size. We will remove this dependence when making a profile decomposition. Then another application of the scaling symmetry $(\ref{1.3})$, this time with $\lambda = 2R$ implies

\begin{equation}\label{3.15}
\| P_{> 2R} u_{0} \|_{\dot{H}^{1/2}(\mathbf{R}^{3})} + \| P_{> 2R} u_{1} \|_{\dot{H}^{-1/2}(\mathbf{R}^{3})} \leq \epsilon,
\end{equation}

\begin{equation}\label{3.16}
\| (1 - \chi(2x)) u_{0} \|_{\dot{H}^{1/2}(\mathbf{R}^{3})} + \| (1 - \chi(2x)) u_{1} \|_{\dot{H}^{-1/2}(\mathbf{R}^{3})} \leq \epsilon,
\end{equation}

\begin{equation}\label{3.17}
\| u \|_{L_{t,x}^{4}([-\frac{\delta}{2R}, \frac{\delta}{2R}] \times \mathbf{R}^{3})} \lesssim \epsilon_{0},
\end{equation}
and finally

\begin{equation}\label{3.18}
\| u \|_{L_{t}^{\infty} B_{1, 2}^{1/2}([-\frac{\delta}{2R}, \frac{\delta}{2R}] \times \mathbf{R}^{3})} \lesssim \| u_{0} \|_{B_{1,1}^{2}(\mathbf{R}^{3})} + \| u_{1} \|_{B_{1,1}^{1}(\mathbf{R}^{3})}.
\end{equation}
The next step is to show that this local solution has a singularity that is isolated in a suitable sense. Observe that the dispersive estimates imply that the linear wave equation $u_{tt} - \Delta u = 0$ with initial data $(u_{0}, u_{1})$ lies in $L^{\infty}$ when $t > 0$. Indeed,

\begin{equation}\label{3.19}
\| S(t)(u_{0}, u_{1}) \|_{L^{\infty}} \lesssim \frac{1}{t} \sum_{j} [2^{2j} \| P_{j} u_{0} \|_{L^{1}(\mathbf{R}^{3})} + 2^{j} \| P_{j} u_{1} \|_{L^{1}(\mathbf{R}^{3})}] \lesssim \frac{1}{t} [\| u_{0} \|_{B_{1,1}^{2}(\mathbf{R}^{3})} + \| u_{1} \|_{B_{1,1}^{1}(\mathbf{R}^{3})}].
\end{equation}
Interpolating $(\ref{3.19})$ with Bernstein's inequality, for any $j$,

\begin{equation}\label{3.20}
\| S(t)(P_{j} u_{0}, P_{j} u_{1}) \|_{L^{6}(\mathbf{R}^{3})} \lesssim \frac{2^{-j/6}}{t^{2/3}} [2^{2j} \| P_{j} u_{0} \|_{L^{1}(\mathbf{R}^{3})} + 2^{j} \| P_{j} u_{1} \|_{L^{1}(\mathbf{R}^{3})}],
\end{equation}
while by the Sobolev embedding theorem $\dot{H}^{1}(\mathbf{R}^{3}) \hookrightarrow L^{6}(\mathbf{R}^{3})$,

\begin{equation}\label{3.21}
\| S(t) P_{j}(u_{0}, u_{1}) \|_{L^{6}(\mathbf{R}^{3})} \lesssim 2^{j/2}  [2^{2j} \| P_{j} u_{0} \|_{L^{1}(\mathbf{R}^{3})} + 2^{j} \| P_{j} u_{1} \|_{L^{1}(\mathbf{R}^{3})}],
\end{equation}
so then by direct computation

\begin{equation}\label{3.22}
\sup_{t > 0} t^{1/2} \| S(t)(u_{0}, u_{1}) \|_{L^{6}(\mathbf{R}^{3})} + \| S(t)(u_{0}, u_{1}) \|_{L_{t}^{2} L_{x}^{6}(\mathbf{R} \times \mathbf{R}^{3})} \lesssim \| u_{0} \|_{B_{1,1}^{2}(\mathbf{R}^{3})} + \| u_{1} \|_{B_{1,1}^{1}(\mathbf{R}^{3})}.
\end{equation}

\begin{lemma}\label{l3.2}
If $\delta > 0$ is given by the local result in lemma $\ref{l3.1}$ for some $\epsilon_{0} > 0$, then

\begin{equation}\label{3.23}
\sup_{-\frac{\delta}{2R} < t < \frac{\delta}{2R}} t^{1/2} \| u \|_{L_{x}^{6}(\mathbf{R}^{3})} + \| u \|_{L_{t}^{2} L_{x}^{6}([-\frac{\delta}{2R}, \frac{\delta}{2R}] \times \mathbf{R}^{3})} \lesssim \| u_{0} \|_{B_{1,1}^{2}(\mathbf{R}^{3})} + \| u_{1} \|_{B_{1,1}^{1}(\mathbf{R}^{3})}.
\end{equation}
\end{lemma}

\noindent \emph{Proof:} By the dispersive estimates (theorem $\ref{t2.2}$), the Hardy - Littlewood - Sobolev inequality, and interpolation

\begin{equation}\label{3.24}
\aligned
\| \int_{0}^{t} S(t - \tau) F(u(\tau)) d\tau \|_{L_{t}^{2} L_{x}^{6}([0, \frac{\delta}{2R}] \times \mathbf{R}^{3})} \lesssim \| |\nabla|^{1/3} F(u) \|_{L_{t,x}^{6/5}([0, \frac{\delta}{2R}] \times \mathbf{R}^{3})} \\
 \lesssim \| |\nabla|^{1/2} u \|_{L_{t}^{\infty} L_{x}^{2}([0, \frac{\delta}{2R}] \times \mathbf{R}^{3})}^{2/3} \| u \|_{L_{t,x}^{4}([0, \frac{\delta}{2R}] \times \mathbf{R}^{3})}^{4/3} \| u \|_{L_{t}^{2} L_{x}^{6}([0, \frac{\delta}{2R}] \times \mathbf{R}^{3})} \lesssim \epsilon_{0}^{4/3} \| u \|_{L_{t}^{2} L_{x}^{6}([0, \frac{\delta}{2R}] \times \mathbf{R}^{3})}.
\endaligned
\end{equation}
Combining $(\ref{3.24})$ with $(\ref{3.22})$ proves

\begin{equation}\label{3.25}
\| u  \|_{L_{t}^{2} L_{x}^{6}([0, \frac{\delta}{2R}] \times \mathbf{R}^{3})} \lesssim \| u_{0} \|_{B_{1,1}^{2}(\mathbf{R}^{3})} + \| u_{1} \|_{B_{1,1}^{1}(\mathbf{R}^{3})}.
\end{equation}
Next let $c > 0$ be a small constant to be determined later. Again by theorem $\ref{t2.2}$, the Hardy - Littlewood - Sobolev inequality, and interpolation,

\begin{equation}\label{3.26}
\aligned
\sup_{t \in [0, \frac{\delta}{2R}]} t^{1/2} \| \int_{0}^{(1 - c)t} S(t - \tau) F(u(\tau)) d\tau \|_{L^{6}(\mathbf{R}^{3})} \lesssim \frac{1}{c^{1/2}} \| |\nabla|^{1/3} F(u) \|_{L_{t,x}^{6/5}([0, \frac{\delta}{2R}] \times \mathbf{R}^{3})} \\ \lesssim \frac{1}{c^{1/2}} \| |\nabla|^{1/2} u \|_{L_{t}^{\infty} L_{x}^{2}([0, \frac{\delta}{2R}] \times \mathbf{R}^{3})}^{2/3} \| u \|_{L_{t,x}^{4}([0, \frac{\delta}{2R}] \times \mathbf{R}^{3})}^{4/3} \| u \|_{L_{t}^{2} L_{x}^{6}([0, \frac{\delta}{2R}] \times \mathbf{R}^{3})} \\ \lesssim \frac{\epsilon_{0}^{4/3}}{c^{1/2}} ( \| u_{0} \|_{B_{1,1}^{2}(\mathbf{R}^{3})} + \| u_{1} \|_{B_{1,1}^{1}(\mathbf{R}^{3})}).
\endaligned
\end{equation}
Also for any $t \in [0, \frac{\delta}{2R}]$, by theorem $\ref{t2.2}$,

\begin{equation}\label{3.27}
\aligned
t^{1/2} \| \int_{(1 - c) t}^{t} S(t - \tau) F(u(\tau)) d\tau \|_{L^{6}(\mathbf{R}^{3})} \lesssim (\sup_{t \in [0, \frac{\delta}{2R}]} t^{1/2} \| u(t) \|_{L^{6}(\mathbf{R}^{3})})^{5/3} \| |\nabla|^{1/2} u \|_{L_{t}^{\infty} L_{x}^{2}([0, \frac{\delta}{2R}] \times \mathbf{R}^{3})}^{2/3} \\ \times \| u \|_{L_{t}^{\infty} L_{x}^{3}([0, \frac{\delta}{2R}] \times \mathbf{R}^{3})}^{2/3} \cdot  \int_{(1 - c)t}^{t} \frac{1}{(t - \tau)^{2/3}} \frac{1}{t^{1/3}} d\tau \lesssim c^{1/3} (\sup_{t \in [0, \frac{\delta}{2R}]} t^{1/2} \| u(t) \|_{L^{6}(\mathbf{R}^{3})})^{5/3}.
\endaligned
\end{equation}
Therefore,

\begin{equation}\label{3.28}
\| u \|_{L_{t}^{2} L_{x}^{6}([0, \frac{\delta}{2R}] \times \mathbf{R}^{3})} + \sup_{0 < t < \delta} t^{1/2} \| u(t) \|_{L^{6}} \lesssim \| u_{0} \|_{B_{1,1}^{2}(\mathbf{R}^{3})} + \| u_{1} \|_{B_{1,1}^{1}(\mathbf{R}^{3})}.
\end{equation}
Then by time reversal symmetry the proof of lemma $\ref{l3.2}$ is complete. $\Box$\vspace{5mm}

\noindent Next, we show that a local solution may be written as a sum of a term with bounded energy and a term with good dispersive properties. To simplify notation let $\delta_{1} = \frac{\delta}{2R}$. By energy inequalities, Strichartz estimates (theorem $\ref{t2.1}$), and lemma $\ref{l3.2}$,

\begin{equation}\label{3.29}
\| \int_{\frac{\delta_{1}}{10}}^{\delta_{1}} S(t - \tau) F(u(\tau)) d\tau \|_{\dot{H}^{1} \times L^{2}(\mathbf{R}^{3})} \lesssim \| u \|_{L_{t}^{3} L_{x}^{6}([\frac{\delta_{1}}{10}, \delta] \times \mathbf{R}^{3})}^{3} \lesssim \frac{1}{\delta_{1}^{1/2}}.
\end{equation}
Next, by the radial Sobolev embedding theorem (theorem $\ref{t2.3}$) and $(\ref{3.13})$, if $\chi \in C_{0}^{\infty}(\mathbf{R}^{3})$ is supported on $|x| \leq 1$, $\chi(x) = 1$ on $|x| \leq \frac{1}{2}$, then

\begin{equation}\label{3.30}
\aligned
\| (1 - \chi(\frac{10x}{\delta_{1}})) F(u) \|_{L_{t}^{1} L_{x}^{2}([0, \frac{\delta_{1}}{10}] \times \mathbf{R}^{3})} \\ \lesssim \delta_{1}^{1/2} \| (1 - \chi(\frac{10 x}{\delta_{1}})) u \|_{L_{t,x}^{\infty}([0, \frac{\delta_{1}}{10}] \times \mathbf{R}^{3})} \| u \|_{L_{t,x}^{4}([0, \frac{\delta_{1}}{10}] \times \mathbf{R}^{3})}^{2} \lesssim \frac{1}{\delta_{1}^{1/2}}.
\endaligned
\end{equation}
Now for $t > \delta_{1}$ let

\begin{equation}\label{3.31}
v(t) = S(t) \chi(\frac{10 x}{\delta})(u_{0}, u_{1}) + \int_{0}^{\delta/10} S(t - \tau) \chi(\frac{10 x}{\delta}) F(u(\tau)) d\tau.
\end{equation}
Combining lemma $\ref{l2.4}$ with

\begin{equation}\label{3.33}
\| [P_{j}, \chi] F(u) \|_{L_{t}^{1} \dot{H}^{-1/2}([-\frac{\delta_{1}}{10}, \frac{\delta_{1}}{10}] \times \mathbf{R}^{3})} \lesssim 2^{-j} \delta_{1}^{-1} \| F(u) \|_{L_{t}^{1} L_{x}^{3/2}([-\frac{\delta_{1}}{10}, \frac{\delta_{1}}{10}] \times \mathbf{R}^{3})},
\end{equation}

\noindent lemma $\ref{l3.2}$,

\begin{equation}\label{3.34}
\| P_{\leq 0} \chi(\frac{10x}{\delta_{1}}) F(u) \|_{\dot{H}^{-1}(\mathbf{R}^{3})} \lesssim \| u \|_{L_{x}^{3}(\mathbf{R}^{3})}^{3},
\end{equation}
$(\ref{3.10})$ - $(\ref{3.13})$, the sharp Huygens principle, which implies $v$ is supported on $\{ (x, t) : ||x| - t| \leq \frac{\delta_{1}}{2} \}$, and the radial Sobolev embedding theorem (theorem $\ref{t2.3}$),

\begin{equation}\label{3.32}
\| v(t) \|_{L^{\infty}(\mathbf{R}^{3})} \lesssim \frac{1}{t} [\| u_{0} \|_{B_{1,1}^{2}(\mathbf{R}^{3})} + \| u_{1} \|_{B_{1,1}^{1}(\mathbf{R}^{3})}].
\end{equation}
This implies good properties of $S(t - \delta_{1})(v(\delta_{1}), v_{t}(\delta_{1}))$.

\begin{lemma}\label{l3.4}
Let $w(\delta_{1}) + v(\delta_{1}) = u(\delta_{1})$. Then

\begin{equation}\label{3.35}
w(\delta_{1}) = S(\delta_{1})(1 - \chi(\frac{10 x}{\delta_{1}}))(u_{0}, u_{1}) + \int_{0}^{\frac{\delta_{1}}{10}} S(\delta_{1} - \tau) (1 - \chi(\frac{10x}{\delta_{1}})) F(u(\tau)) d\tau + \int_{\frac{\delta_{1}}{10}}^{\delta_{1}} S(\delta_{1} - \tau) F(u(\tau)) d\tau,
\end{equation}
and

\begin{equation}\label{3.36}
\| w(\delta_{1}) \|_{\dot{H}^{1} \times L^{2}(\mathbf{R}^{3})} \lesssim \delta_{1}^{-1/2}.
\end{equation}
\end{lemma}
\emph{Proof:} By $(\ref{3.30})$ and $(\ref{3.31})$ it only remains to compute

\begin{equation}\label{3.37}
\| (1 - \chi(\frac{10x}{\delta_{1}})) u_{0} \|_{\dot{H}^{1}(\mathbf{R}^{3})} + \| (1 - \chi(\frac{10x}{\delta_{1}})) u_{1} \|_{L^{2}(\mathbf{R}^{3})}.
\end{equation}
First,

\begin{equation}\label{3.38}
|u_{1}(0, r)|  \lesssim \int_{r}^{\infty} |\partial_{r} u_{1}(0, s)| ds \lesssim \frac{1}{r^{2}},
\end{equation}
so

\begin{equation}\label{3.39}
\int_{\frac{\delta_{1}}{10}}^{\infty} |u_{1}(r, 0)|^{2} r^{2} dr \lesssim \int_{\frac{\delta_{1}}{10}}^{\infty} |u_{1}(r, 0)| dr \lesssim \frac{1}{\delta_{1}}.
\end{equation}
Next, for any $j$, by the Sobolev embedding theorem,

\begin{equation}\label{3.40}
\int_{0}^{2^{-j}} \frac{r^{2}}{r} | \partial_{r} (P_{j} u_{0})| dr \lesssim 2^{2j} \| P_{j} u_{0} \|_{L^{1}(\mathbf{R}^{3})},
\end{equation}
while by Bernstein's inequality

\begin{equation}\label{3.41}
\int_{2^{-j}}^{\infty} \frac{r^{2}}{r} |\partial_{r} (P_{j} u_{0})| dr \lesssim 2^{2j} \| P_{j} u_{0} \|_{L^{1}(\mathbf{R}^{3})}.
\end{equation}
Therefore, $\| \frac{1}{r} \partial_{r} u_{0} \|_{L^{1}(\mathbf{R}^{3})} \lesssim \| u_{0} \|_{B_{1,1}^{2}(\mathbf{R}^{3})}$, and since $u_{0}$ is radially symmetric $\Delta u_{0} = (\partial_{rr} + \frac{2}{r} \partial_{r}) u_{0}$, so $\| \partial_{rr} u_{0} \|_{L^{1}(\mathbf{R}^{3})} \lesssim \| u_{0} \|_{B_{1,1}^{2}(\mathbf{R}^{3})}$. By the fundamental theorem of calculus,

\begin{equation}\label{3.42}
|u_{r}(0, r)| \leq \int_{r}^{\infty} |u_{rr}(0, s)| ds \lesssim \frac{1}{r^{2}}.
\end{equation}
Therefore,

\begin{equation}\label{3.43}
\int_{\frac{\delta_{1}}{10}}^{\infty} |u_{r}(0, r)|^{2} r^{2} dr \lesssim \int_{\frac{\delta_{1}}{10}}^{\infty} |u_{r}(0, r)| dr \lesssim \frac{1}{\delta_{1}},
\end{equation}
and

\begin{equation}\label{3.44}
\int_{0}^{\infty} (u_{r}(r))^{2} r^{3} dr \leq \int_{0}^{\infty} (\int_{r}^{\infty} |u_{rr}(s)| ds) r dr \lesssim \int_{0}^{\infty} |u_{rr}(s)| s^{2} ds < \infty.
\end{equation}
$\Box$

\section{Proof of global well - posedness}
In this section we extend local well - posedness to global well - posedness, proving

\begin{theorem}\label{t4.1}
$(\ref{1.1})$ is globally well - posed, and for any compact interval $J \subset \mathbf{R}$,

\begin{equation}\label{4.0}
\| u \|_{L_{t,x}^{4}(J \times \mathbf{R}^{3})} < \infty.
\end{equation}
\end{theorem}
\emph{Proof:} By time reversal symmetry, to prove this it suffices to show that the local well - posedness result of lemma $\ref{l3.1}$ can be extended to all times $t > \delta_{1}$. Throughout the proof the implicit constant depends on $\delta_{1}$ and $\| u_{0} \|_{B_{1,1}^{2}} + \| u_{1} \|_{B_{1,1}^{1}}$. Now for $t > \delta_{1}$ let

\begin{equation}\label{4.1}
u(t) = w(t) + v(t),
\end{equation}
where $v(t)$ is given by $(\ref{3.31})$ and $w$ solves

\begin{equation}\label{4.2}
w_{tt} - \Delta w = -u^{3}.
\end{equation}


\noindent Next, copying $(\ref{1.2})$ let $E(w(t))$ be the energy of $w$,

\begin{equation}\label{4.3}
E(w(t)) = \frac{1}{2} \int |\nabla w(t,x)|^{2} dx + \frac{1}{2} \int (w_{t}(t,x))^{2} dx + \frac{1}{4} \int (w(t,x))^{4} dx.
\end{equation}
By $(\ref{3.18})$, $(\ref{3.36})$, and the Sobolev embedding theorem $w \in L^{3} \cap L^{6}$, so 

\begin{equation}\label{4.4}
E(w(\delta_{1})) \lesssim 1.
\end{equation}
Next,

\begin{equation}\label{4.5}
\aligned
\frac{d}{dt} E(w(t)) = \int ((w(t,x))^{3} - (u(t,x))^{3}) w_{t}(t,x) dx \\ = - \int w_{t}(t,x) [(v(t,x))^{3} + 3 v(t,x)^{2} w(t,x) + 3 v(t,x) w(t,x)^{2}] dx.
\endaligned
\end{equation}
Now by $(\ref{3.32})$,

\begin{equation}\label{4.6}
\int w_{t}(t,x) (w(t,x))^{2} v(t,x) dx \lesssim \| v(t) \|_{L^{\infty}(\mathbf{R}^{3})} \| w(t) \|_{L^{4}(\mathbf{R}^{3})}^{2} \| w_{t}(t) \|_{L^{2}(\mathbf{R}^{3})} \lesssim \frac{1}{t} E(w(t)).
\end{equation}

\begin{equation}\label{4.7}
\int w_{t}(t,x) (v(t,x))^{3} dx \lesssim \| w_{t}(t) \|_{L^{2}(\mathbf{R}^{3})} \| v(t) \|_{L^{\infty}(\mathbf{R}^{3})} \| v(t) \|_{L^{4}(\mathbf{R}^{3})}^{2} \lesssim \frac{1}{t} E(w(t))^{1/2} \| v(t) \|_{L^{4}(\mathbf{R}^{3})}^{2}.
\end{equation}
Finally

\begin{equation}\label{4.8}
\aligned
\int w_{t}(t,x) v(t,x)^{2} w(t,x) dx \lesssim \| w_{t}(t) \|_{L^{2}(\mathbf{R}^{3})} \| w(t) \|_{L^{4}(\mathbf{R}^{3})} \| v(t) \|_{L^{4}(\mathbf{R}^{3})} \| v(t) \|_{L^{\infty}(\mathbf{R}^{3})} \\
\lesssim \frac{1}{t} E(w(t))^{3/4} \| v(t) \|_{L^{4}(\mathbf{R}^{3})}.
\endaligned
\end{equation}
Then by interpolation

\begin{equation}\label{4.9}
\frac{d}{dt} E(w(t)) \lesssim \frac{1}{t} E(w(t)) + \frac{1}{t} \| v(t) \|_{L^{4}(\mathbf{R}^{3})}^{4}.
\end{equation}
By Strichartz estimates (theorem $\ref{t2.1}$) $v \in L_{t,x}^{4}(\mathbf{R} \times \mathbf{R}^{3})$. Then by Gronwall's inequality and time reversal symmetry there exist constants $C_{1}(\| u_{0} \|_{B_{1,1}^{2}}, \| u_{1} \|_{B_{1,1}^{1}}, \delta_{1})$ and $C_{2}(\| u_{0} \|_{B_{1,1}^{2}}, \| u_{1} \|_{B_{1,1}^{1}}, \delta_{1})$ such that

\begin{equation}\label{4.10}
E(w(t)) \lesssim C_{1} (1 + |t|)^{C_{2}}.
\end{equation}
Notice that this implies that for any compact interval $J \subset \mathbf{R}$,

\begin{equation}\label{4.11}
\| w \|_{L_{t,x}^{4}(J \times \mathbf{R}^{3})} + \| v \|_{L_{t,x}^{4}(J \times \mathbf{R}^{3})} < \infty.
\end{equation}
This proves the theorem. $\Box$

\section{Hyperbolic coordinates}
In this section we prove

\begin{theorem}\label{t5.1}
The global solution given in theorem $\ref{t4.1}$ scatters both forward and backward in time.
\end{theorem}
\noindent \emph{Proof:} By time reversal symmetry and $(\ref{1.13})$, it suffices to show that

\begin{equation}\label{5.0}
\int_{0}^{\infty} \int_{0}^{\infty} u(t,r)^{4} r^{2} dr dt < \infty.
\end{equation}
First we make a translation in time so that $t_{0} = 0$ maps to $t_{0} = 1 - \delta_{1}$. Theorem $\ref{t4.1}$ implies that

\begin{equation}\label{5.0.1}
\int_{0}^{1} \int_{0}^{\infty} u(t,r)^{4} r^{2} dr dt < \infty.
\end{equation}
Next, by small data arguments (see for example \cite{LS}), the solution to $(\ref{1.1})$ with initial data given by $(\ref{3.14})$, has finite $L_{t,x}^{4}$ norm. Then by finite propagation speed this implies

\begin{equation}\label{5.0.2}
\int_{1}^{\infty} \int_{r > \frac{1}{2} + t} u(t,r)^{4} r^{2} dr dt \lesssim \epsilon_{0}.
\end{equation}
It therefore remains to prove

\begin{equation}\label{5.0.3}
\int_{1}^{\infty} \int_{r \leq \frac{1}{2} + t} u(t,r)^{4} r^{2} dr dt < \infty.
\end{equation}
It is convenient to compute this norm in hyperbolic coordinates.

\begin{equation}\label{5.1}
\tilde{w}(\tau, s) = \frac{e^{\tau} \sinh s}{s} w(e^{\tau} \cosh s, e^{\tau} \sinh s),
\end{equation}

\begin{equation}\label{5.2}
\tilde{v}(\tau, s) = \frac{e^{\tau} \sinh s}{s} v(e^{\tau} \cosh s, e^{\tau} \sinh s),
\end{equation}
and

\begin{equation}\label{5.3}
\tilde{u}(\tau, s) = \frac{e^{\tau} \sinh s}{s} u(e^{\tau} \cosh s, e^{\tau} \sinh s).
\end{equation}
Then $\tilde{w}$ solves the nonlinear wave equation

\begin{equation}\label{5.4}
\partial_{\tau \tau} \tilde{w} - \partial_{ss} \tilde{w} - \frac{2}{s} \partial_{s} \tilde{w} = -(\frac{s}{\sinh s})^{2} \tilde{u}^{3}.
\end{equation}
We begin by showing that $\tilde{w}$ has finite hyperbolic energy.

\begin{lemma}\label{l5.1}
There exists some $0 < \tau < \delta_{1}$ such that
\begin{equation}\label{5.5}
E(\tilde{w}(\tau)) < \infty,
\end{equation}
where $E(\tilde{w}(\tau))$ is the hyperbolic energy

\begin{equation}\label{5.6}
\frac{1}{2} \int_{0}^{\infty} \tilde{w}_{s}(\tau, s)^{2} s^{2} ds + \frac{1}{2} \int_{0}^{\infty} \tilde{w}_{\tau}(\tau, s)^{2} s^{2} ds + \frac{1}{4} \int_{0}^{\infty} (\frac{s}{\sinh s})^{2} \tilde{w}(\tau, s)^{4} s^{2} ds.
\end{equation}
\end{lemma}
\emph{Proof:} When $\tau = 0$,

\begin{equation}\label{5.7}
\lim_{s \rightarrow \infty} e^{\tau} \cosh s - e^{\tau} \sinh s = 1.
\end{equation}
Combining

\begin{equation}\label{5.8}
v(t) = S(t) \chi(\frac{10 x}{\delta_{1}})(u_{0}, u_{1}) + \int_{0}^{\delta/10} S(t - \tau) \chi(\frac{10 x}{\delta_{1}}) F(u(\tau)) d\tau
\end{equation}
with finite propagation speed, there exists some $s_{0}$ such that

\begin{equation}\label{5.9}
\int_{s_{0}}^{\infty} s^{2} \tilde{w}_{s}(\tau, s)^{2} ds + \int_{s_{0}}^{\infty} s^{2} \tilde{w}_{\tau}(\tau, s) ds = \int_{s_{0}}^{\infty} s^{2} \tilde{u}_{s}(\tau, s)^{2} ds + \int_{s_{0}}^{\infty} s^{2} \tilde{u}_{\tau}(\tau, s) ds.
\end{equation}
Now by standard properties of the wave equation, remembering that $t_{0} = 1 - \delta_{1}$,

\begin{equation}\label{5.10}
\aligned
s \tilde{u}(\tau, s) = \frac{1}{2} (e^{\tau + s} - (1 - \delta_{1})) u_{0}(e^{\tau + s} - (1 - \delta_{1})) + \frac{1}{2} (1 - \delta_{1} - e^{\tau - s}) u_{0}(1 - \delta_{1} - e^{\tau - s}) \\+ \frac{1}{2} \int_{e^{\tau - s} + (1 - \delta_{1})}^{e^{\tau + s} - (1 - \delta_{1})} u_{1}(r) r dr
+ \frac{1}{2} \int_{1 - \delta_{1}}^{e^{\tau} \cosh s} \int_{-e^{\tau - s} + t}^{e^{\tau + s} - t} r u^{3}(t, r) dr dt.
\endaligned
\end{equation}
Then

\begin{equation}\label{5.11}
\aligned
\partial_{\tau + s} (s \tilde{u}(\tau,s))|_{\tau = 0} = \frac{1}{2} e^{s} u_{0}(e^{s} - (1 - \delta_{1})) \\ + \frac{1}{2} e^{s} (e^{s} - (1 - \delta_{1})) u_{0}'(e^{s} - (1 - \delta_{1})) \\
+ \frac{e^{s}}{2} (e^{s} - (1 - \delta_{1})) u_{1}(e^{s} - (1 - \delta_{1})) \\ + \frac{e^{s}}{2} \int_{1 - \delta_{1}}^{\cosh s} (e^{s} - t) u^{3}(t, e^{s} - t) dt,
\endaligned
\end{equation}
and

\begin{equation}\label{5.11}
\aligned
\partial_{\tau - s} (s \tilde{u}(\tau,s))|_{\tau = 0} = \frac{1}{2} e^{-s} u_{0}((1 - \delta_{1}) - e^{-s}) \\ + \frac{e^{-s}}{2} ((1 - \delta_{1}) - e^{-s}) u_{0}'((1 - \delta_{1}) - e^{-s}) \\
+ \frac{e^{-s}}{2} ((1 - \delta_{1}) - e^{-s}) u_{1}((1 - \delta_{1}) - e^{-s}) \\ + \frac{e^{-s}}{2} \int_{1 - \delta_{1}}^{e^{\tau} \cosh s} (t - e^{-s}) u^{3}(t, t - e^{-s}) dt.
\endaligned
\end{equation}
First, making a change of variables and using $(\ref{3.38})$ - $(\ref{3.44})$,

\begin{equation}\label{5.12}
\int_{s_{0}}^{\infty} e^{2s} u_{0}(e^{s} - (1 - \delta_{1}))^{2} ds \lesssim \int_{0}^{\infty} u_{0}(r)^{2} r dr < \infty,
\end{equation}

\begin{equation}\label{5.13}
\int_{s_{0}}^{\infty} e^{2s} (e^{s} - (1 - \delta_{1}))^{2} u_{0}'(e^{s} - (1 - \delta_{1}))^{2} ds \lesssim \int_{0}^{\infty} (\partial_{r} u_{0}(r))^{2} r^{3} dr < \infty,
\end{equation}
and

\begin{equation}\label{5.14}
\int_{s_{0}}^{\infty} e^{2s} (e^{s} - (1 - \delta_{1}))^{2} u_{1}(e^{s} - (1 - \delta_{1}))^{2} ds \lesssim \int_{0}^{\infty} r^{3} u_{1}(r)^{2} dr < \infty.
\end{equation}

\noindent Also, by $(\ref{3.38})$ - $(\ref{3.44})$ combined with the fact that $\cosh s - \sinh s \geq \frac{1}{2}$ when $s \geq s_{0}$, $(\ref{5.0.1})$, and $(\ref{5.0.2})$,

\begin{equation}\label{5.15}
\int_{s_{0}}^{\infty} e^{2s} (\int_{1 - \delta_{1}}^{\cosh s} (e^{s} - t) u^{3}(t, e^{s} - t) dt)^{2} ds \lesssim \int_{s_{0}}^{\infty} \int_{1 - \delta_{1}}^{\cosh s} e^{3s} (e^{s} - t)^{2} u^{6}(t, e^{s} - t) dt ds < \infty.
\end{equation}

\begin{equation}\label{5.16}
\int_{s_{0}}^{\infty} e^{-2s} u_{0}((1 - \delta_{1}) - e^{-s})^{2} ds \lesssim \int_{s_{0}}^{\infty} e^{-2s} ds < \infty.
\end{equation}

\begin{equation}\label{5.17}
\int_{s_{0}}^{\infty} e^{-2s} ((1 - \delta_{1}) - e^{-s})^{2} (u_{0}'((1 - \delta_{1}) - e^{-s}))^{2} ds \lesssim \int_{s_{0}}^{\infty} e^{-2s} ds < \infty.
\end{equation}

\begin{equation}\label{5.18}
\int_{s_{0}}^{\infty} e^{-2s} ((1 - \delta_{1}) - e^{-s})^{2} u_{1}((1 - \delta_{1}) - e^{-s})^{2} ds \lesssim \int_{s_{0}}^{\infty} e^{-2s} ds < \infty.
\end{equation}
Also by $(\ref{3.38})$ - $(\ref{3.44})$,

\begin{equation}\label{5.19}
\int_{s_{0}}^{\infty} e^{-2s} ( \int_{1 - \delta_{1}}^{e^{\tau} \cosh s} (t - e^{-s}) u^{3}(t, t - e^{-s}) dt)^{2} ds \lesssim \int_{s_{0}}^{\infty} e^{-2s} ds < \infty.
\end{equation}
In fact the above computations could be made for any $0 < \tau < \delta_{1}$ with some uniform $s_{0}$. So to prove the lemma it suffices to show that

\begin{equation}\label{5.19.1}
\int_{0}^{\delta_{1}} \int_{0}^{s_{0}} s^{2} \tilde{w}_{s}(\tau, s)^{2} ds d\tau + \int_{0}^{\delta_{1}} \int_{0}^{s_{0}} s^{2} \tilde{w}_{\tau}(\tau, s) ds d\tau.
\end{equation}
This fact is an immediate consequence of $(\ref{5.1})$, theorem $\ref{t4.1}$, and the fact that $e^{\tau} \sinh s$ and $e^{\tau} \cosh s$ are uniformly bounded when $s \leq s_{0}$ and $\tau \leq \delta_{1}$. Thus, for some $0 < \tau_{0} < \delta_{1}$,

\begin{equation}\label{5.19.2}
\int_{0}^{\infty} s^{2} \tilde{w}_{s}(\tau_{0}, s)^{2} ds + \int_{0}^{\infty} s^{2} \tilde{w}_{\tau}(\tau_{0}, s) ds < \infty.
\end{equation}
Then an application of the Sobolev embedding theorem completes the proof of lemma $\ref{l5.1}$. $\Box$\vspace{5mm}

\noindent Next we compute

\begin{equation}\label{5.20}
\frac{d}{d\tau} E(\tilde{w}(\tau)) = \int \tilde{w}_{\tau} [\tilde{u}^{3} - \tilde{w}^{3}] (\frac{s}{\sinh s})^{2} s^{2} ds.
\end{equation}
By the support properties of $v$ and the definition of $\tilde{v}$,

\begin{equation}\label{5.21}
\| \tilde{w}_{\tau}(\tau) \|_{L^{2}} \| \tilde{v}(\tau, s)^{2} (\frac{s}{\sinh s}) \|_{L^{2}} \| \tilde{v}(\tau,s) (\frac{s}{\sinh s}) \|_{L^{\infty}} \lesssim e^{-\tau/2} E(\tilde{w}(\tau))^{1/2} \| \tilde{v}(\tau, s)^{2} (\frac{s}{\sinh s}) \|_{L^{2}}.
\end{equation}
Meanwhile,

\begin{equation}\label{5.22}
\| \tilde{w}_{\tau}(\tau) \|_{L^{2}} \| \tilde{v}(\tau,s) (\frac{s}{\sinh s}) \|_{L^{\infty}} \| \tilde{w}(\tau,s)^{2} (\frac{s}{\sinh s}) \|_{L^{2}} \lesssim E(\tilde{w}(\tau)) e^{-\tau/2},
\end{equation}
and

\begin{equation}\label{5.23}
\aligned
\| \tilde{w}_{\tau}(\tau) \|_{L^{2}} \| \tilde{v}(\tau,s) (\frac{s}{\sinh s}) \|_{L^{\infty}} \| \tilde{w}(\tau,s) (\frac{s}{\sinh s})^{1/2} \|_{L^{4}} \| \tilde{v}(\tau,s) (\frac{s}{\sinh s})^{1/2} \|_{L^{4}} \\ \lesssim e^{-\tau/2} E(\tilde{w}(\tau))^{3/4} \| \tilde{v}(\tau,s) (\frac{s}{\sinh s})^{1/2} \|_{L^{4}}.
\endaligned
\end{equation}
Now by a change of variables

\begin{equation}\label{5.24}
\aligned
\int_{0}^{\infty} \int_{0}^{\infty} \tilde{v}(\tau, s)^{4} s^{2} (\frac{s}{\sinh s})^{2} ds d\tau \\
= \int_{0}^{\infty} \int_{0}^{\infty} \frac{e^{4 \tau} (\sinh s)^{4}}{s^{4}} (\frac{s^{4}}{(\sinh s)^{2}}) v(e^{\tau} \cosh s, e^{\tau} \sinh s)^{4} ds d\tau \\
= \int \int_{t^{2} - r^{2} \geq 1} v(t,r)^{4} r^{2} dr dt.
\endaligned
\end{equation}
Therefore, since $\| v \|_{L_{t,x}^{4}} < \infty$, Gronwall's inequality implies that $E(\tilde{w}(\tau))$ is uniformly bounded on $\mathbf{R}$.\vspace{5mm}

\noindent Next we prove the Morawetz estimate. 

\begin{theorem}\label{t5.2}
\begin{equation}\label{5.25}
\int \int \tilde{w}(s, \tau)^{4} (\frac{s}{\sinh s})^{2} s^{2} ds d\tau < \infty.
\end{equation}
\end{theorem}
\noindent \emph{Proof:} Let

\begin{equation}\label{5.26}
M(\tau) = \int \tilde{w}_{\tau} (\frac{x}{|x|} \cdot \nabla \tilde{w}) dx.
\end{equation}
Then

\begin{equation}\label{5.27}
\frac{d}{d\tau} M(\tau) = \int (\frac{\cosh s}{\sinh s}) (\frac{s}{\sinh s})^{2} \tilde{w}^{4} s^{2} ds + \int \frac{x}{|x|} \cdot (\nabla \tilde{w}) (\tilde{u}^{3} - \tilde{w}^{3}) s^{2} ds d\tau.
\end{equation}
As in the bounded energy computations,

\begin{equation}\label{5.28}
\| \tilde{w}_{s} \|_{L^{2}} \| \tilde{v}^{2} (\frac{s}{\sinh s}) \|_{L^{2}} \| \tilde{v} (\frac{s}{\sinh s}) \|_{L^{\infty}} \lesssim e^{-\tau/2} E(w(\tau))^{1/2} \| \tilde{v}^{2} (\frac{s}{\sinh s}) \|_{L^{2}},
\end{equation}

\begin{equation}\label{5.29}
\| w_{s} \|_{L^{\infty}} \| \tilde{v} (\frac{s}{\sinh s}) \|_{L^{\infty}} \| \tilde{w}^{2} (\frac{s}{\sinh s}) \|_{L^{2}} \lesssim E(\tilde{w}) e^{-\tau/2},
\end{equation}
and

\begin{equation}\label{5.30}
\aligned
\| w_{s} \|_{L^{\infty}} \| \tilde{v} (\frac{s}{\sinh s}) \|_{L^{\infty}} \| \tilde{w} (\frac{s}{\sinh s})^{1/2} \|_{L^{4}} \| \tilde{v} (\frac{s}{\sinh s})^{1/2} \|_{L^{4}} \\ \lesssim E(\tilde{w})^{3/4} e^{-\tau/2} \| \tilde{v} (\frac{s}{\sinh s})^{1/2} \|_{L^{4}}.
\endaligned
\end{equation}
Therefore, by the fundamental theorem of calculus, the fact that the energy is uniformly bounded, and $(\ref{5.27})$,

\begin{equation}\label{5.31}
\int_{0}^{\infty} \int_{0}^{\infty} \tilde{w}(s, \tau)^{4} (\frac{s}{\sinh s})^{2} s^{2} ds d\tau < \infty.
\end{equation}
\noindent $\Box$\vspace{5mm}

\noindent Then by the change of variables in $(\ref{5.24})$, $(\ref{5.31})$ implies

\begin{equation}\label{5.32}
\int \int_{t^{2} - r^{2} \geq 1} w(t,r)^{4} r^{2} dr dt < \infty.
\end{equation}
Finally, by theorem $\ref{t4.1}$,

\begin{equation}\label{5.33}
\int_{1}^{\infty} \int_{t^{2} - r^{2} \geq 1, t - r \leq \frac{1}{2}} u(t,r)^{4} r^{2} dr dt < \infty.
\end{equation}
Combining $(\ref{5.0.1})$, $(\ref{5.0.2})$, $(\ref{5.24})$, $(\ref{5.32})$, and $(\ref{5.33})$,

\begin{equation}\label{5.34} 
\int_{0}^{\infty} \int_{0}^{\infty} u(t,r)^{4} r^{2} dr dt < \infty.
\end{equation}
This proves theorem $\ref{t5.1}$. $\Box$\vspace{5mm}

\noindent \textbf{Remark:} Notice that theorem $\ref{t5.1}$ implies that

\begin{equation}\label{5.35}
\int_{0}^{\infty} \int_{0}^{\infty} u(t,r)^{4} r^{2} dr dt \leq C(\| u_{0} \|_{B_{1,1}^{2}}, \| u_{1} \|_{B_{1,1}^{1}}, \delta_{1}) < \infty.
\end{equation}
Thus theorem $\ref{t5.1}$ is not equivalent to theorem $\ref{t1.6}$. This $\delta_{1} > 0$ depends on the support of $u_{0}$ and $u_{1}$ in space $(\ref{3.14})$ and in frequency $(\ref{3.2})$. To remove this requirement, it is necessary make a profile decomposition, the subject of the final section of this paper.

\section{Profile decomposition}
Observe that the difficulty in going theorem $\ref{t5.1}$ to theorem $\ref{t1.6}$ lies in the fact that even if most of the $B_{1,1}^{2} \times B_{1,1}^{1}$ norm lies below frequency one, guaranteeing local well - posedness on an interval of length $2 \delta$ (lemma $\ref{l3.1}$), the $R$ appearing in $(\ref{3.14})$ could be very large, and thus after rescaling, $\delta_{1} > 0$ could be quite small.\vspace{5mm}

\noindent However, the intuition guiding an important refinement utilizes finite propagation speed. Indeed, if $(u_{0}, u_{1})$ were radial functions supported on the annulus $R \leq r \leq 2R$, $R$ large, then the $\dot{H}^{1/2}$ norm on balls of radius $c R$ for some $c > 0$ small would actually be fairly small. Therefore, one could then apply the small data arguments of \cite{LS} to prove the result of lemma $\ref{l3.1}$ actually holds on an interval of length $\sim R$.\vspace{5mm}

\noindent By the uncertainty principle, when most of the $B_{1,1}^{2} \times B_{1,1}^{1}$ lies below frequency one we have $R \gtrsim 1$. Utilizing a by now standard profile decomposition argument, it is possible to show that the above argument can refine theorem $\ref{t5.1}$ when $R >> 1$, thus proving theorem $\ref{t1.6}$. A key ingredient is the profile decomposition of \cite{Ramos}.

\begin{theorem}[Profile decomposition]\label{t6.1}
Suppose that there is a uniformly bounded, radially symmetric sequence

\begin{equation}\label{6.1}
\| u_{0}^{n} \|_{\dot{H}^{1/2}(\mathbf{R}^{3})} + \| u_{1}^{n} \|_{\dot{H}^{-1/2}(\mathbf{R}^{3})} \leq C_{0} < \infty.
\end{equation}
Then there exists a subsequence, also denoted $(u_{0}^{n}, u_{1}^{n}) \subset \dot{H}^{1/2} \times \dot{H}^{-1/2}$ such that for any $N < \infty$,

\begin{equation}\label{6.2}
S(t)(u_{0}^{n}, u_{1}^{n}) = \sum_{j = 1}^{N} \Gamma_{j}^{n} S(t)(\phi_{0}^{j}, \phi_{1}^{j}) + S(t)(R_{0, n}^{N}, R_{1,n}^{N}),
\end{equation}
with

\begin{equation}\label{6.3}
\lim_{N \rightarrow \infty} \limsup_{n \rightarrow \infty} \| S(t)(R_{0,n}^{N}, R_{1,n}^{N}) \|_{L_{t,x}^{4}(\mathbf{R} \times \mathbf{R}^{3})} = 0.
\end{equation}
$\Gamma_{j}^{n}$ is the action of the group $(0, \infty) \times \mathbf{R}$,

\begin{equation}\label{6.4}
\Gamma_{j}^{n} F(t,x) = \lambda_{j}^{n} F(\lambda_{j}^{n} (t - t_{n}^{j}), \lambda_{j}^{n} x).
\end{equation}
Additionally, for every $j \neq k$,

\begin{equation}\label{6.5}
\lim_{n \rightarrow \infty} \frac{\lambda_{j}^{n}}{\lambda_{k}^{n}} + \frac{\lambda_{k}^{n}}{\lambda_{j}^{n}} + (\lambda_{n}^{j})^{1/2} (\lambda_{n}^{k})^{1/2} |t_{j}^{n} - t_{k}^{n}| = \infty.
\end{equation}
Furthermore, for every $N \geq 1$,

\begin{equation}\label{6.6}
\| (u_{0, n}, u_{1, n}) \|_{\dot{H}^{1/2} \times \dot{H}^{-1/2}}^{2} = \sum_{j = 1}^{N} \| (\phi_{0}^{j}, \phi_{0}^{k}) \|_{\dot{H}^{1/2} \times \dot{H}^{-1/2}}^{2} + \| (R_{0, n}^{N}, R_{1, n}^{N}) \|_{\dot{H}^{1/2} \times \dot{H}^{-1/2}}^{2} + o_{n}(1).
\end{equation}
\end{theorem}
\textbf{Remark:} \cite{Ramos} proved this result for data which need not be radially symmetric. Such a result is substantially more difficult since it requires accounting for Lorentz transformations and translation in space. See \cite{Gerard} and \cite{BG} for the early development of the profile decomposition.\vspace{5mm}

\noindent Now let

\begin{equation}
\aligned
f(M) = \sup \{ \| u \|_{L_{t,x}^{4}(\mathbf{R} \times \mathbf{R}^{3})} : u \text{ solves } (\ref{1.1}) \text{ with initial data } \\ (u_{0}, u_{1}) \in B_{1,1}^{2} \times B_{1,1}^{1}, \hspace{5mm} \| u_{0} \|_{B_{1,1}^{2}} + \| u_{1} \|_{B_{1,1}^{1}} \leq M \}.
\endaligned
\end{equation}
To prove theorem $\ref{t1.6}$ it suffices to show that $f(M) < \infty$ for any $M$.\vspace{5mm}

\noindent \textbf{Remark:} Theorem $\ref{t1.5}$ implies that such a function is well - defined.\vspace{5mm}

\noindent Take a uniformly bounded sequence

\begin{equation}\label{6.7}
\| u_{0}^{n} \|_{B_{1,1}^{2}(\mathbf{R}^{3})} + \| u_{1}^{n} \|_{B_{1,1}^{1}(\mathbf{R}^{3})} \leq C_{0} < \infty,
\end{equation}
such that if $u^{n}(t)$ is the solution of $(\ref{1.1})$ with initial data $(u_{0}^{n}, u_{1}^{n})$, then 

\begin{equation}
\| u^{n}(t) \|_{L_{t,x}^{4}(\mathbf{R} \times \mathbf{R}^{3})} \rightarrow f(C_{0}).
\end{equation}
By the Sobolev embedding theorem

\begin{equation}\label{6.8}
\| u_{0}^{n} \|_{\dot{H}^{1/2}(\mathbf{R}^{3})} + \| u_{1}^{n} \|_{\dot{H}^{-1/2}(\mathbf{R}^{3})} \lesssim \| u_{0}^{n} \|_{B_{1,1}^{2}(\mathbf{R}^{3})} + \| u_{1}^{n} \|_{B_{1,1}^{1}(\mathbf{R}^{3})} \leq C_{0} < \infty,
\end{equation}
which by theorem $\ref{t6.1}$ gives a profile decomposition

\begin{equation}\label{6.9}
S(t)(u_{0}^{n}, u_{1}^{n}) = \sum_{j = 1}^{N} S(t - t_{n}^{j}) (\lambda_{n}^{j} \phi_{0}^{j}(\lambda_{n}^{j} x), (\lambda_{n}^{j})^{2} \phi_{1}^{j}(\lambda_{n}^{j} x)) + S(t)(R_{0, n}^{N}, R_{1,n}^{N}).
\end{equation}
In the course of proving theorem $\ref{t6.1}$, \cite{Ramos} proved

\begin{equation}\label{6.10}
S(-\frac{t_{n}^{j}}{\lambda_{n}^{j}})(\frac{1}{\lambda_{n}^{j}} u_{0}^{n}(\frac{x}{\lambda_{n}^{j}}), \frac{1}{(\lambda_{n}^{j})^{2}} u_{1}^{n}(\frac{x}{\lambda_{n}^{j}})) \rightharpoonup \phi_{0}^{j}(x)
\end{equation}
weakly in $\dot{H}^{1/2}(\mathbf{R}^{3})$, and

\begin{equation}\label{6.10}
\partial_{t}S(t -\frac{t_{n}^{j}}{\lambda_{n}^{j}})(\frac{1}{\lambda_{n}^{j}} u_{0}^{n}(\frac{x}{\lambda_{n}^{j}}), \frac{1}{(\lambda_{n}^{j})^{2}} u_{1}^{n}(\frac{x}{\lambda_{n}^{j}}))|_{t = 0} \rightharpoonup \phi_{0}^{j}(x)
\end{equation}
weakly in $\dot{H}^{-1/2}(\mathbf{R}^{3})$.

\begin{lemma}\label{l6.2}
For each $j$, $\frac{t_{n}^{j}}{\lambda_{n}^{j}}$ is uniformly bounded.
\end{lemma}
\emph{Proof:} By the dispersive estimates (theorem $\ref{t2.2}$), for any $l \in \mathbf{Z}$,

\begin{equation}\label{6.11}
\aligned
\| P_{l} S(t - \frac{t_{n}^{j}}{\lambda_{n}^{j}}) ((\lambda_{n}^{j})^{-1} u_{0}^{n}(\frac{x}{\lambda_{n}^{j}}), (\lambda_{n}^{j})^{-2} u_{1}^{n}(\frac{x}{\lambda_{n}^{j}})) \|_{L^{\infty}(\mathbf{R}^{3})} \\ \lesssim \frac{1}{|t - \frac{t_{n}^{j}}{\lambda_{n}^{j}}|} [2^{2l} \| P_{l} ((\lambda_{n}^{j})^{-1} u_{0}^{n}(\frac{x}{\lambda_{n}^{j}}) \|_{L^{1}(\mathbf{R}^{3})} + 2^{l} \| P_{l} ((\lambda_{n}^{j})^{-2} u_{1}^{n}(\frac{x}{\lambda_{n}^{j}})) \|_{L^{1}(\mathbf{R}^{3})}].
\endaligned
\end{equation}
Meanwhile, by Bernstein's inequality and the Sobolev embedding theorem

\begin{equation}\label{6.12}
\aligned
\| P_{l} S(t - \frac{t_{n}^{j}}{\lambda_{n}^{j}}) ((\lambda_{n}^{j})^{-1} u_{0}^{n}(\frac{x}{\lambda_{n}^{j}}), (\lambda_{n}^{j})^{-2} u_{1}^{n}(\frac{x}{\lambda_{n}^{j}})) \|_{L^{2}(\mathbf{R}^{3})} \\ \lesssim 2^{-l/2} [2^{2l} \| P_{l} ((\lambda_{n}^{j})^{-1} u_{0}^{n}(\frac{x}{\lambda_{n}^{j}}) \|_{L^{1}(\mathbf{R}^{3})} + 2^{l} \| P_{l} ((\lambda_{n}^{j})^{-2} u_{1}^{n}(\frac{x}{\lambda_{n}^{j}})) \|_{L^{1}(\mathbf{R}^{3})}].
\endaligned
\end{equation}
Then by interpolation, for any $l \in \mathbf{Z}$,

\begin{equation}\label{6.13}
\aligned
\| P_{l} S(t - \frac{t_{n}^{j}}{\lambda_{n}^{j}}) ((\lambda_{n}^{j})^{-1} u_{0}^{n}(\frac{x}{\lambda_{n}^{j}}), (\lambda_{n}^{j})^{-2} u_{1}^{n}(\frac{x}{\lambda_{n}^{j}})) \|_{L_{t,x}^{4}(\{ |t - \frac{t_{n}^{j}}{\lambda_{n}^{j}}| > C 2^{-l} \} \times \mathbf{R}^{3})} \\ \lesssim \frac{1}{C^{1/4}} [2^{2l} \| P_{l} ((\lambda_{n}^{j})^{-1} u_{0}^{n}(\frac{x}{\lambda_{n}^{j}}) \|_{L^{1}(\mathbf{R}^{3})} + 2^{l} \| P_{l} ((\lambda_{n}^{j})^{-2} u_{1}^{n}(\frac{x}{\lambda_{n}^{j}})) \|_{L^{1}(\mathbf{R}^{3})}].
\endaligned
\end{equation}
Then if $\limsup_{n \rightarrow \infty} \frac{|t_{n}^{j}|}{\lambda_{n}^{j}} = \infty$, then possibly after passing to a subsequence,

\begin{equation}\label{6.14}
S(t - \frac{t_{n}^{j}}{\lambda_{n}^{j}}) ((\lambda_{n}^{j})^{-1} u_{0}^{n}(\frac{x}{\lambda_{n}^{j}}), (\lambda_{n}^{j})^{-2} u_{1}^{n}(\frac{x}{\lambda_{n}^{j}})) \rightharpoonup 0
\end{equation}
weakly in $L_{t,x}^{4}(\mathbf{R} \times \mathbf{R}^{3})$. Utilizing lemma $4.1$ from \cite{Ramos},

\begin{lemma}\label{l6.3}
\begin{equation}\label{6.15}
(u_{0, n}, u_{1,n}) \rightharpoonup (\phi_{0}, \phi_{1})
\end{equation}
weakly in $\dot{H}^{1/2}(\mathbf{R}^{3}) \times \dot{H}^{-1/2}(\mathbf{R}^{3})$ is equivalent to

\begin{equation}\label{6.16}
S(t)(u_{0, n}, u_{1,n}) \rightharpoonup S(t)(\phi_{0}, \phi_{1})
\end{equation}
weakly in $L_{t,x}^{4}(\mathbf{R} \times \mathbf{R}^{3})$.
\end{lemma}
Therefore $(\ref{6.14})$ implies $(\phi_{0}^{j}, \phi_{1}^{j}) = (0, 0)$. $\Box$\vspace{5mm}

\noindent Therefore $\frac{|t_{n}^{j}|}{\lambda_{n}^{j}}$ is uniformly bounded, so after passing to a subsequence, $\frac{t_{n}^{j}}{\lambda_{n}^{j}}$ converges to some $t_{0}^{j} \in \mathbf{R}$. Then

\begin{equation}\label{6.17}
S(\frac{t_{n}^{j}}{\lambda_{n}^{j}})(\phi_{0}^{j}, \phi_{1}^{j}) \rightarrow S(t_{0}^{j})(\phi_{0}^{j}, \phi_{1}^{j})
\end{equation}
strongly in $\dot{H}^{1/2}(\mathbf{R}^{3}) \times \dot{H}^{-1/2}(\mathbf{R}^{3})$. Absorbing the error into $(R_{0, n}^{N}, R_{1,n}^{N})$ and taking

\begin{equation}\label{6.18}
(\tilde{\phi}_{0}^{j}, \tilde{\phi}_{1}^{j}) = S(t_{0}^{j})(\phi_{0}, \phi_{1}),
\end{equation}
we can assume $t_{n}^{j} \equiv 0$. Therefore,

\begin{equation}\label{6.19}
(u_{0}^{n}, u_{1}^{n}) = \sum_{j = 1}^{N} (\lambda_{n}^{j} \phi_{0}^{j}(\lambda_{n}^{j} x), (\lambda_{n}^{j})^{2} \phi_{1}^{j}(\lambda_{n}^{j} x)) + (R_{0,n}^{N}, R_{1,n}^{N}),
\end{equation}
and

\begin{equation}\label{6.20}
\lim_{n \rightarrow \infty} \frac{\lambda_{j}^{n}}{\lambda_{k}^{n}} + \frac{\lambda_{k}^{n}}{\lambda_{j}^{n}} = \infty.
\end{equation}
But then

\begin{equation}\label{6.21}
\| u_{0}^{n} \|_{B_{1,1}^{2}(\mathbf{R}^{3})} + \| u_{1}^{n} \|_{B_{1,1}^{1}(\mathbf{R}^{3})} \leq C_{0} < \infty
\end{equation}
combined with lemma $\ref{l6.3}$, $(\ref{6.19})$, and $(\ref{6.20})$ implies that for any $j$,

\begin{equation}\label{6.22}
\| \phi_{0}^{j} \|_{B_{1,1}^{2}(\mathbf{R}^{3})} + \| \phi_{1}^{j} \|_{B_{1,1}^{1}(\mathbf{R}^{3})} \leq C_{0}.
\end{equation}
Possibly reordering $j$, $(\ref{6.6})$ implies that there exists $N_{0}(\epsilon, C_{0})$ such that if $j \geq N_{0}(\epsilon)$,

\begin{equation}\label{6.23}
\| (\phi_{0}^{j}, \phi_{1}^{j}) \|_{\dot{H}^{1/2} \times \dot{H}^{-1/2}} < \epsilon.
\end{equation}
Now for each $j$ let $v^{j}(t,x)$ be the solution of $(\ref{1.1})$ with initial data $(\phi_{0}^{j}, \phi_{1}^{j})$. By the small data arguments of \cite{LS}, when $j \geq N_{0}(\epsilon)$,

\begin{equation}\label{6.24}
\| v^{j} \|_{L_{t,x}^{4}(\mathbf{R} \times \mathbf{R}^{3})} \lesssim \| \phi_{0}^{j} \|_{\dot{H}^{1/2}(\mathbf{R}^{3})} + \| \phi_{1}^{j} \|_{\dot{H}^{-1/2}(\mathbf{R}^{3})}.
\end{equation}
Meanwhile, by theorem $\ref{t5.1}$ combined with $(\ref{6.22})$, when $j \leq N_{0}(\epsilon)$,

\begin{equation}\label{6.25}
\| v^{j} \|_{L_{t,x}^{4}(\mathbf{R} \times \mathbf{R}^{3})} \lesssim_{j, C_{0}} 1.
\end{equation}
Also by $(\ref{6.20})$, for any $j \neq k$, the Lebesgue dominated convergence theorem implies

\begin{equation}\label{6.26}
\lim_{n \rightarrow \infty} \int \int |\lambda_{n}^{j} v^{j}(\lambda_{n}^{j} t, \lambda_{n}^{j} x)|^{2} |\lambda_{n}^{k} v^{k}(\lambda_{n}^{k} t, \lambda_{n}^{k} x)|^{2} dx dt = 0.
\end{equation}
Therefore,

\begin{equation}\label{6.27}
\lim_{n \rightarrow \infty} \| \sum_{1 \leq j \leq N} \lambda_{n}^{j} v^{j}(\lambda_{n}^{j} t, \lambda_{n}^{j} x) \|_{L_{t,x}^{4}(\mathbf{R} \times \mathbf{R}^{3})}
\end{equation}
is uniformly bounded, independent of $N$. Also,

\begin{equation}\label{6.28}
\aligned
F(\sum_{j = 1}^{N} \lambda_{n}^{j} v^{j}(\lambda_{n}^{j} t, \lambda_{n}^{j} x)) - \sum_{j = 1}^{N} F(\lambda_{n}^{j} v^{j}(\lambda_{n}^{j} t, \lambda_{n}^{j} x)) \\ = \sum_{1 \leq j \neq k \leq N} O(|\lambda_{n}^{j} v^{j}(\lambda_{n}^{j} t, \lambda_{n}^{j} x)| |\lambda_{n}^{k} v^{k}(\lambda_{n}^{k} t, \lambda_{n}^{k} x)|^{2}),
\endaligned
\end{equation}
so by $(\ref{6.25})$, $(\ref{6.26})$, and $(\ref{6.27})$,

\begin{equation}\label{6.29}
\lim_{n \rightarrow \infty} \| F(\sum_{j = 1}^{N} \lambda_{n}^{j} v^{j}(\lambda_{n}^{j} t, \lambda_{n}^{j} x)) - \sum_{j = 1}^{N} F(\lambda_{n}^{j} v^{j}(\lambda_{n}^{j} t, \lambda_{n}^{j} x)) \|_{L_{t,x}^{4/3}(\mathbf{R} \times \mathbf{R}^{3})} = 0.
\end{equation}
Therefore, by lemma $\ref{l2.5}$, the solution $u^{n}_{N}(t,x)$ to $(\ref{1.1})$ with initial data

\begin{equation}\label{6.30}
\sum_{j = 1}^{N} (\lambda_{n}^{j} \phi_{0}^{j}(\lambda_{n}^{j} x), (\lambda_{n}^{j})^{2} \phi_{1}^{j}(\lambda_{n}^{j} x))
\end{equation}
has

\begin{equation}\label{6.31}
\lim_{n \rightarrow \infty} \| u_{N}^{n}(t) \|_{L_{t,x}^{4}(\mathbf{R} \times \mathbf{R}^{3})}
\end{equation}
bounded uniformly in $N$. By another application of lemma $\ref{l2.5}$ combined with

\begin{equation}\label{6.32}
\lim_{N \rightarrow \infty} \limsup_{n \rightarrow \infty} \| S(t)(R_{0, n}^{N}, R_{1,n}^{N}) \|_{L_{t,x}^{4}(\mathbf{R} \times \mathbf{R}^{3})} = 0,
\end{equation}
if $u^{n}(t)$ is the solution to $(\ref{1.1})$ with initial data $(u_{0}^{n}, u_{1}^{n})$ satisfying $(\ref{6.1})$, then

\begin{equation}\label{6.33}
\| u^{n}(t) \|_{L_{t,x}^{4}(\mathbf{R} \times \mathbf{R}^{3})}
\end{equation}
is uniformly bounded. This proves theorem $\ref{t1.6}$. $\Box$

\nocite*
\bibliographystyle{plain}

\end{document}